\documentclass{amsart}
\usepackage{amssymb}

\newcommand{\Z}{{\mathbb{Z}}}
\newcommand{\Q}{{\mathbb{Q}}}
\newcommand{\N}{{\mathbb{N}}}

\newcommand{\bA}{{\mathbf{A}}}
\newcommand{\bC}{{\mathbf{C}}}
\newcommand{\bH}{{\mathbf{H}}}
\newcommand{\bM}{{\mathbf{M}}}
\newcommand{\bP}{{\mathbf{P}}}
\newcommand{\bR}{{\mathbf{R}}}
\newcommand{\bT}{{\mathbf{T}}}
\newcommand{\fC}{{\mathfrak{C}}}
\newcommand{\cE}{{\mathcal{E}}}
\newcommand{\cL}{{\mathcal{L}}}

\renewcommand{\leq}{\leqslant}
\renewcommand{\geq}{\geqslant}
\renewcommand{\atop}[2]{\genfrac{}{}{0pt}{}{#1}{#2}}

\newtheorem{thm}{Theorem}[section]
\newtheorem{lem}[thm]{Lemma}
\newtheorem{cor}[thm]{Corollary}
\newtheorem{prop}[thm]{Proposition}
\newtheorem{conj}[thm]{Conjecture}

\theoremstyle{definition}
\newtheorem{exmp}[thm]{Example}
\newtheorem{defn}[thm]{Definition}

\theoremstyle{remark}
\newtheorem{rem}[thm]{Remark}

\begin{document}

\title{Computing Kazhdan--Lusztig cells for unequal parameters}

\author{Meinolf Geck}
\address{Institut Girard Desargues, bat. Jean Braconnier, Universit\'e Lyon 1, 
21 av Claude Bernard, F--69622 Villeurbanne Cedex, France}
\email{geck@desargues.univ-lyon1.fr}

\date{November 2003}
\begin{abstract} Following Lusztig, we consider a Coxeter group $W$ together 
with a weight function $L$. This gives rise to the pre-order relation 
$\leq_{L}$ and the corresponding partition of $W$ into left cells. We 
introduce an equivalence relation on weight functions such that, in 
particular, $\leq_{L}$ is constant on equivalent classes. We shall work 
this out explicitly for $W$ of type $F_4$ and check that several of Lusztig's 
conjectures concerning left cells with unequal parameters hold in this case,
even for those parameters which do not admit a geometric interpretation.
The proofs involve some explicit computations using {\sf CHEVIE}.
\end{abstract}

\maketitle

\section{Introduction} \label{sec:intro}
This paper is concerned with the computation of the Kazhdan--Lusztig 
polynomials, the left cells and the corresponding representations of 
a finite Coxeter group $W$ with respect to a weight function $L$. 
Following Lusztig \cite{Lusztig03}, a weight function on $W$ is a function
$L\colon W\rightarrow\Z$ such that $L(ww')=L(w)+L(w')$ whenever $l(ww')=
l(w)+l(w')$ where $l$ is the length function on $W$. As in most parts of
\cite{Lusztig03}, we shall only consider weight functions such that 
$L(w)>0$ for all $w\neq 1$. 

The case where $L$ is constant on the generators of $W$ is known as
the equal parameter case. If, moreover, $W$ is a finite Weyl group,
then there is a geometric interpretation for the Kazhdan--Lusztig
polynomials and this leads to many deep properties for which no
elementary proofs are known (see \cite{Lu1}, \cite{Lu2}). Recently, Lusztig 
\cite{Lusztig03} has formulated a number of precise conjectures
in the general case of unequal parameters. Furthermore, Lusztig proposes
a geometric interpretation at least for those weight functions which
arise in the representation theory of finite groups of Lie type. (The
complete list of these $L$ is given in \cite{Lusztig77}, Table~II,  p.~35.) 

One of our aims here is to show that some of Lusztig's conjectures hold
for $W$ of type $F_4$ and any weight function, even for those $L$ which
do not admit a geometric interpretation. In type $F_4$, with generators
and diagram given by the diagram below, a weight function $L$ is specified
by two positive integers $a:=L(s_1)=L(s_2)>0$ and $b:=L(s_3)=L(s_4)>0$. 
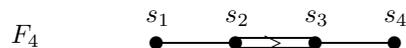
\begin{center}
\begin{picture}(200,20)
\put( 10, 5){$F_4$} 
\put( 61,13){$s_1$}
\put( 91,13){$s_2$}
\put(121,13){$s_3$}
\put(151,13){$s_4$}
\put( 65, 5){\circle*{5}}
\put( 95, 5){\circle*{5}}
\put(125, 5){\circle*{5}}
\put(155, 5){\circle*{5}}
\put(105,2.5){$>$}
\put( 65, 5){\line(1,0){30}}
\put( 95, 7){\line(1,0){30}}
\put( 95, 3){\line(1,0){30}}
\put(125, 5){\line(1,0){30}}
\end{picture}
\end{center}
By explicit computations  using the {\sf CHEVIE}-system, we obtain the
following results.

\begin{thm} \label{typf4a} Let $W$ be of type $F_4$ and $L$ any weight 
function on $W$ with $L(w)>0$ for $w\neq 1$. Then the left cell
representations of $W$ (with respect to $L$) are precisely the 
constructible representations, as defined by Lusztig 
\cite[Chap.~22]{Lusztig03}.  
\end{thm}

The above result is conjectured to hold in general by Lusztig 
\cite[\S 22.29]{Lusztig03}. As  far as the partition of $W$ into left cells
is concerned, we shall see that there are only four essentially different 
cases, according to whether $b=a$, $b=2a$, $2a>b>a$ or $b>2a$; see
Corollary~\ref{corconf4} and Remark~\ref{rem0}.

\begin{thm} \label{typf4b} Let be of type $F_4$ and $L$ any weight 
function on $W$ with $L(w)>0$ for $w\neq 1$. For $w\in W$, we define
$\Delta(w)\in \Z_{\geq 0}$ and $0 \neq n_w\in \Z$ by the condition 
\[ P_{1,w}^*=n_wv^{-\Delta(w)}+\mbox{strictly smaller powers of $v$};
\quad \text{see Lusztig \protect{\cite[14.1]{Lusztig03}}}.\]
Let $C$ be a left cell $C$ of $W$ (with respect to $L$). Then the function 
$w\mapsto \Delta(w)$ reaches its minumum at exactly one element of $C$,
denoted by $d_C\in C$. We have $d_C^{\, 2}=1$ and $n_{d_C}=\pm 1$.
\end{thm}

(For the definition of $P_{y,w}^*$, see Section~2.) The elements $d_C$  
are the {\em distinguished involutions} whose existence is predicted by
Lusztig \cite{Lusztig03}, Conjectures~14.2 (P1, P6, P13). The following 
result is also part of those conjectures (P4, P9).

\begin{thm} \label{typf4c} Let be of type $F_4$ and $L$ any weight 
function on $W$ with $L(w)>0$ for $w\neq 1$. For any $y,w\in W$, we have
the following implication:
\[ y \leq_L w \qquad \mbox{and}\qquad y \sim_{LR} w \qquad 
\Longrightarrow \qquad y\sim_L w.\]
\end{thm}

(For the definition of the relations $\leq_L$, $\sim_L$, $\sim_{LR}$, 
see Section~2.) The proofs of the above three theorems will be given 
in Section~4 (see Corollary~\ref{corconf4}).

In type $F_4$, there is a geometric interpretation for the cases where 
$(a,b)\in \{(1,1), (1,2),(1,4)\}$; see \cite{Lusztig77}, Table~II (p.~35). 
To deal with arbibrary values for $a$ and $b$, we have to provide a 
theoretical argument which shows that it is enough to consider only those 
$L$ where the values on the generators are bounded by a constant which can 
be explicitly computed in terms of $W$. More precisely, in 
Definition~\ref{def1}, we introduce (for general $W$) an equivalence 
relation on the set of weight functions. This relation has the property 
that two related weight functions give rise to the same partition of $W$ 
into left cells, the same left pre-order relation and the same set of left 
cell representations. In Corollary~\ref{plem3}, we show that any weight 
function is equivalent to a weight function whose values on the generators 
are bounded by a constant which can be computed efficiently. 

Lusztig's results \cite{Lusztig03} on dihedral groups are interpreted
in this framework in Example~\ref{exp1}. Conjecture~\ref{conj1} (found 
independently by Bonnaf\'e) would yield a complete description of the 
equivalence classes of weight functions in type~$B_n$. 

Both the results in type $F_4$ and the evidence for the conjecture on type
$B_n$ are based on a {\sf CHEVIE}-program which we have developed, for 
computing the Kazhdan--Lusztig polynomials, the $M$-polynomials, and the 
pre-order relations $\leq_L$, $\leq_{LR}$ for a finite Coxeter group $W$ 
and any choice of the parameters (either given by independent indeterminates 
and a monomial order on them, or given by a weight function). For example, 
this program systematically computes the polynomials $P_{y,w}^*$ for all 
pairs $y<w$ in $W$; it also computes all incidences of the Kazhdan--Lusztig
pre-order relation $y\leq_L w$.  The program automatically checks some of 
Lusztig's conjectures (in particular, the properties expressed in the 
above three theorems) and computes the characters carried by the various 
left cells. 
These programs have already been used
in the computations reported in \cite[\S 11.3]{ourbuch} and 
\cite[\S 7]{myert02}. To my knowledge, the first such programs 
(for Kazhdan--Lusztig polynomials in the unequal parameter case) were 
written by K.~Bremke \cite{Br} who used them to compute $W$-graphs for the 
irreducible representations of certain Iwahori--Hecke algebras of 
type $F_4$. We only remark that, in the case of equal parameters, there
is already a rather sophisticated theory for the computation of 
Kazhdan--Lusztig polynomials; see Alvis \cite{Al} and Ducloux \cite{Ducloux}.

\section{Total orderings and weight functions} \label{sec:left}
The basic references for this section are \cite{Lusztig83} and 
\cite{Lusztig03}. In the latter reference, Lusztig studies the left
cells of a Coxeter group $W$ with respect to a weight function $L$ on
$W$. In the former reference, Lusztig considers a more abstract setting
where left cells are defined with respect to an abelian group and a total 
order on it. We will see in this section that the more abstract setting can 
be used to show that two given weight functions actually give rise to 
the same partition of $W$ into left cells. (A similar argument has already 
been used, for example, in \cite{BI}.) This will provide the theoretical 
argument for showing that, in order to determine the left cells for all 
possible weight functions on $W$, it is actually enough to consider a 
certain finite number of weight functions.

We begin by recalling the basic setting for the definition of 
Kazhdan--Lusztig polynomials and left cells. Let $W$ be a Coxeter group, 
with generating set $S$. Let $\Gamma$ be an abelian group (written 
multiplicatively) and $\bA={\Z}[\Gamma]$ be the group algebra of $\Gamma$ 
over $\Z$. Let $\{v_s\mid s\in S\}\subset \Gamma$ be a subset such
that $v_s=v_t$ whenever $s,t\in S$ are conjugate in $W$. Then we have a
corresponding generic Iwahori--Hecke algebra $\bH$, with $\bA$-basis 
$\{\bT_w\mid w \in W\}$ and multiplication given by the rule
\begin{equation*} \bT_s\,\bT_w=\left\{\begin{array}{cl} \bT_{sw} & 
\quad \mbox{if $l(sw)>l(w)$},\\
\bT_{sw}+(v_s-v_s^{-1})\bT_w & \quad 
\mbox{if $l(sw)<l(w)$};\end{array} \right. \tag{2.1}
\end{equation*} 
here, $l\colon W \rightarrow {\N}_0$ denotes the usual length function
on $W$ with respect to $S$. (Note that the above elements $\bT_w$ are
denoted $\tilde{T}_w$ in \cite{Lusztig83}.) 

Let $a \mapsto \bar{a}$ be the involution of ${\Z}[\Gamma]$ which takes
$g$ to $g^{-1}$ for any $g \in \Gamma$. We extend it to a map $\bH\rightarrow 
\bH$, $h \mapsto \overline{h}$, by the formula
\begin{equation*}
\overline{\sum_{w \in W} a_w \bT_w}=\sum_{w \in W} \bar{a}_w
\bT_{w^{-1}}^{-1} \qquad (a_w \in {\Z}[\Gamma]).\tag{2.2}
\end{equation*}
Then $h \mapsto \overline{h}$ is in fact a ring involution.  

Now assume that we have chosen a total ordering of $\Gamma$. This is
specified by a multiplicatively closed subset $\Gamma_{+} \subseteq \Gamma 
\setminus \{1\}$ such that we have $\Gamma=\Gamma_{+} \amalg \{1\} \amalg 
\Gamma_{-}$, where $\Gamma_{-}= \{g^{-1} \mid g \in \Gamma_{+}\}$. 
Furthermore, we assume that
\begin{equation*}
\{v_s\mid s\in S\} \subset \Gamma_+.\tag{2.3}
\end{equation*}
Given a total ordering of $\Gamma$ as above, we have a corresponding 
{\em Kazhdan--Lusztig basis} of $\bH$, which we denote by $\{\bC_w\mid 
w\in W\}$. (Note that this basis is denoted by $C_w'$ in \cite{Lusztig83}.)
The basis element $\bC_w$ is uniquely determined by the conditions that 
\begin{equation*}
\overline{\bC}_w=\bC_w\qquad \mbox{and}\qquad \bC_w=\bT_w+
\sum_{\atop{y \in W}{y<w}} \bP^{*}_{y,w}\, \bT_y\tag{2.4}
\end{equation*} 
where $\bP^{*}_{y,w}\in {\Z}[\Gamma_-]$ for $y<w$. Here, 
$\leq$ denotes the Bruhat--Chevalley order on $W$. We shall also set 
$\bP_{w,w}^*=1$ for all $w\in W$.  For any $w\in W$ we set $v_w:= v_{s_1}
\cdots v_{s_p}$ where $w=s_1 \cdots s_p$ with $s_i \in S$ is a reduced 
expression. Then we actually have:
\begin{equation*}
\bP_{y,w}:=v_wv_y^{-1}\, \bP_{y,w}^* \mbox{ lies in ${\Z}[v_t^2 \mid 
t\in S]$ and has constant term~$1$}; \tag{2.5}
\end{equation*}
see Lemma~\ref{lem3} below. We have the following multiplication 
formulas. Let $w\in W$ and $s\in S$. Then 
\begin{equation*}
 \bT_s\bC_w= \left\{\begin{array}{cl} \displaystyle \bC_{sw}-v_s\bC_w+
\sum_{\atop{y<w}{sy<y}} \bM_{y,w}^s \, \bC_z &\qquad \mbox{if $sw>w$},\\
v_s\bC_w &\qquad \mbox{if $sw<w$}, \end{array}\right.\tag{2.6}
\end{equation*}
where the coefficients $\bM_{y,w}^s\in \bA$ are such that 
$\overline{\bM}_{y,w}^s=\bM_{y,w}^s$.
Given $y,w\in W$ and $s\in W$, we write $y\leq_{L,s} w$ if the following
conditions are satisfied:
\begin{equation*}
 w=sy>y \quad \mbox{or} \quad sy<y<w<sw \mbox{ and } \bM_{y,w}^s \neq 0.
\tag{2.7}
\end{equation*}
The Kazhdan--Lusztig left preorder $\leq_L$ is the transitive closure
of the above relation, that is, given $y,w\in W$ we have $y\leq_L w$ if
$y=w$ or if there exists a sequence $y=y_0,y_1,\ldots,y_n=w$ of elements in
$W$ and a sequence $s_1,\ldots,s_n$ of generators in $S$ such that
$y_{i-1}\leq_{L,s_i} y_i$ for $1\leq i \leq n$. (See \cite[\S 6]{Lusztig83}.)
Thus, we have $\bH\,\bC_w\subseteq \sum_{y\leq_L w} \bA \bC_y$ for any 
$w\in W$.  The equivalence relation associated with $\leq_L$ will be denoted 
by $\sim_L$ and the corresponding equivalence classes are called the 
{\em left cells} of $W$. Similarly, we write $y\leq_{LR} w$ if $y=w$ or if 
there is a chain of elements $y=y_0,y_1,\ldots,y_n=w$ in $W$ such that, for
each $i$, we have $y_{i-1}\leq_{L} y_i$ or $y_{i-1}^{-1}\leq_L y_i^{-1}$.
The equivalence relation associated with $\leq_{LR}$ will be denoted 
by $\sim_{LR}$ and the corresponding equivalence classes are called the 
{\em two-sided cells} of $W$.  Each two-sided cell  is a union of left
cells and a union of right cells. Consider the following statement:
\begin{equation*}
y \leq_L w \quad \mbox{and} \quad y\sim_{LR} w \quad 
\Longrightarrow \quad y\sim_L w.  \tag{L}
\end{equation*}
This is known to be true in those cases where there is a geometric 
interpretation for the parameters (for example, the equal-parameter
case where $v_s=v_t$ for all $s\neq t$ in $S$); see 
\cite[Chap.~14]{Lusztig03} for more details. The above property plays
an important role in certain representation-theoretic constructions;
see \cite[Chap.~5]{LuBook}. Lusztig \cite[14.2]{Lusztig03} conjectures 
that (L) holds in the general unequal parameter case. It would imply
that the two-sided cells are the {\em minimal} subsets of $W$ which are at
the same time unions of left cells and union of right cells.

Each left cell $\fC$ gives rise to a representation
of~$\bH$. This is constructed as follows (see \cite[\S 7]{Lusztig83}). 
Let $V_{\fC}$ be an $\bA$-module with a free $\bA$-basis $\{e_w \mid w 
\in \fC\}$. Then the action of $\bT_s$ ($s \in S$) is given by the formula
\begin{equation*}
\bT_s.e_w = \left\{\begin{array}{ll} \displaystyle{e_{sw}+v_s e_w-
\sum_{\atop{z<w}{sz<z}} (-1)^{l(w)-l(z)} \bM_{z,w}^s e_z} &\quad \mbox{if 
$sw>w$},\\-v_s^{-1}\,e_w &\quad \mbox{if $sw<w$},\end{array}\right.\tag{2.8}
\end{equation*}
where we tacitly assume that $e_{z}=0$ if $z\not\in \fC$.  (The formula 
(2.8) can be related to the formula (2.6) using a suitable automorphism 
of $\bH$; see \cite[\S 6]{Lusztig83}.) Assume now that $W$ is finite. Upon 
specialization $v_s \mapsto 1$ ($s \in S$), we obtain a representation of 
$W$ which is called the representation carried by $\fC$.  We denote by
$\chi_{\fC}$ the character of that representation, that is, the map 
$w\mapsto \mbox{trace}(w|V_{\fC})$. On the other hand, let $\mbox{Con}
(W,\Gamma_+)$ be the set of so-called constructible characters of $W$, as 
defined by Lustig; see \cite[Chap.~22]{Lusztig03} (and also 
\cite[\S 3]{myert02}, for the general setting with respect to $\Gamma_+
\subset \Gamma$). Consider the following statement:
\begin{equation*}
\mbox{Con}(W,\Gamma_+)=\{\chi_{\fC} \mid \mbox{$\fC$ left cell in $W$ 
with respect to $\Gamma_+\subset \Gamma$}\}.\tag{C}
\end{equation*}
It is conjectured by Lusztig \cite[22.29]{Lusztig03} that (C) always holds. 
This is known to be true in the equal parameter case (see \cite{Lusztig86}) 
and some cases with unequal parameters.  (See, for example, the explicit 
results on type $I_2(m)$ in \cite{Lusztig03}, on type $B_n$ in \cite{BI}, 
and on type $F_4$ in \cite{myert02}). The important point about (C) is that
the constructible characters can be easily determined by a recursive
procedure, using the induction of characters from parabolic subgroups of~$W$.

{\bf Summary.} Given an abelian group $\Gamma$ with a total order specified 
by $\Gamma_{+} \subset \Gamma$ and a choice of parameters $\{v_s\mid s\in S\}
\subset \Gamma_+$, we obtain
\begin{itemize}
\item a collection of polynomials $\bP_{y,w}^* \in {\Z}[\Gamma_-]$ for 
all $y<w$ in $W$;
\item a collection of polynomials $\bM_{y,w}^s \in {\Z}[\Gamma]$ whenever 
$sy<y<w<sw$.
\end{itemize}
These data determine, in a purely combinatorial way, a pre-order relation 
$\leq_L$ on $W$ and the corresponding partition of $W$ into left cells 
and two-sided cells. Finally, we obtain a set of characters of $W$ (the 
characters carried by the left cells).

Now let us specialise the above setting to the case where the parameters
of the Iwahori--Hecke algebra are given by a weight function. Following 
\cite{Lusztig03}, a weight function on $W$ is a function $L\colon W 
\rightarrow \Z$ such that $L(ww')=L(w)+L(w')$ for all $w,w'\in W$ such that 
$l(ww')=l(w)+l(w')$. Such a function is determined by its values $L(s)$ on 
$S$ which are subject only to the condition that $L(s)=L(s')$ for any 
$s\neq s'$ in $S$ such that the order of $ss'$ is finite and odd. (See
Matsumoto's Lemma \cite[\S 1.2]{ourbuch}.) We shall only consider weight 
functions $L$ such that $L(s)>0$ for all $s\in S$. Let $A={\Z}[v,v^{-1}]$ 
where $v$ is an indeterminate. We have a corresponding Iwahori--Hecke 
algebra $H$ with parameters $\{v^{L(s)} \mid s\in S\}$.   Thus, $H$ has
an $A$-basis $\{T_w\mid w\in W\}$ and the multiplication is determined by
the formula
\begin{equation*} 
T_s\,T_w=\left\{\begin{array}{cl} T_{sw} & \quad \mbox{if $l(sw)>l(w)$},\\
T_{sw}+\big(v^{L(s)}-v^{-L(s)}\big)T_w & \quad \mbox{if $l(sw)<l(w)$};
\end{array} \right. \tag{2.9}
\end{equation*}
Now consider the abelian group $\{v^n \mid n \in \Z\}$ with the total order 
specified by $\{v^n\mid n>0\}$. Thus, as above, we have a corresponding
Kazhdan--Lusztig basis $\{C_w \mid w\in W\}$ of $H$. Consequently, we obtain 
\begin{itemize}
\item a collection of polynomials $P_{y,w}^* \in v^{-1}{\Z}[v^{-1}]$ for 
all $y<w$ in $W$;
\item a collection of polynomials $M_{y,w}^s \in {\Z}[v,v^{-1}]$ whenever 
$sy<y<w<sw$.
\end{itemize}
As before, these data determine a pre-order relation $\leq_L$ on $W$ and 
the corresponding partition of $W$ into left cells and two-sided cells;
furthermore, we obtain the characters carried by the left cells of $W$.

The following result establishes a link between the above two situations, 
where we have an abelian group $\Gamma$ with a total order specified by 
$\Gamma_+ \subset \Gamma$ and a choice of parameters $\{v_s\mid s\in S\}
\subset \Gamma_+$ on the one hand, and a weight function $L$ on the other 
hand. As above, denote by $\bP_{y,w}^*$ and $\bM_{y,w}^s$ the
polynomials in ${\Z}[\Gamma]$ arising in the first case, and denote by 
$P_{y,w}^*$ and $M_{y,w}^s$ the polynomials in ${\Z}[v,v^{-1}]$ arising
in the second case. 

We now define two  subsets $\Gamma_+^{(a)}(W),\Gamma_+^{(b)}\subseteq 
\Gamma_+$. First, let $\Gamma_+^{(a)}(W)$ be the of all $\gamma\in 
\Gamma_+$ such that $\gamma^{-1}$ occurs with non-zero coefficient in a
polynomial $\bP_{y,w}^*$ for some $y<w$ in $W$. Next, for any $y,w\in W$ 
and $s\in S$ such that $\bM_{y,w}^s\neq 0$, we write $\bM_{y,w}^s=
n_1\gamma_1+\cdots + n_r\gamma_r$ where $0\neq n_i\in \Z$, $\gamma_i\in 
\Gamma$ and $\gamma_{i-1}^{-1}\gamma_i\in \Gamma_+$ for $2\leq i \leq r$. 
We let $\Gamma_+^{(b)}(W)$ be the set of all elements $\gamma_{i-1}^{-1}
\gamma_i\in \Gamma_+$ arising in this way, for any $y,w,s$ such that 
$\bM_{y,w}^s\neq 0$. Finally, we set $\Gamma_+(W):=\Gamma_+^{(a)}(W)
\cup \Gamma_+^{(b)}(W)$.

\addtocounter{thm}{9}
\begin{prop} \label{prop1} Assume that we have a ring homomorphism
\[ \sigma \colon {\Z}[\Gamma] \rightarrow {\Z}[v,v^{-1}], \qquad v_s 
\mapsto v^{L(s)} \;\;(s\in S)\]
such that 
\begin{equation*}
\sigma(\Gamma_+(W)) \subseteq \{v^n\mid n>0\}. \tag{$*$}
\end{equation*}
Then $\sigma(\bP_{y,w}^*)=P_{y,w}^*$ for all $y<w$ in $W$ and 
$\sigma(\bM_{y,w}^s)=M_{y,w}^s$ for any $s\in S$ such that $sy<y<w<sw$. 
Furthermore, the relations $\leq_L$, $\sim_L$, $\leq_{LR}$
and $\leq_{LR}$ on $W$ defined with respect to the weight function $L$ are 
the same as those with respect to $\Gamma_+\subset \Gamma$, and so are the 
corresponding representations of $W$.
\end{prop}

\begin{proof} The map $\sigma$ induces a ring homomorphism 
\[ \hat{\sigma} \colon \bH \rightarrow H, \qquad \sum_w a_w \,\bT_w \mapsto 
\sum_w {\sigma}(a_w)\,T_w.\]
We have $\overline{\hat{\sigma}(h)}=\hat{\sigma}
(\overline{h})$ for all $h\in \bH$. Thus, applying $\hat{\sigma}$ to 
(2.4),  we obtain 
\[ \overline{\hat{\sigma}(\bC_w)}=\hat{\sigma}(\bC_w) \qquad 
\mbox{and}\qquad \hat{\sigma}(\bC_w)= T_w+ \sum_{\atop{y \in W}{y<w}} 
{\sigma}(\bP^{*}_{y,w})\, T_y\] 
for any $w\in W$. Now condition ($*$) implies that $\sigma(\Gamma_-)
\subseteq \{v^n \mid n<0\}$ and so ${\sigma}(\bP_{y,w}^*)$ is
either $0$ or an integral linear combination of terms $v^n$ with $n<0$.
Thus, the elements $\hat{\sigma}(\bC_w)$ satisfy the defining properties 
for the Kazhdan--Lusztig basis of $H$ and so we must have $\hat{\sigma}
(\bC_w)=C_w$ for all $w\in W$. This also shows that ${\sigma}
(\bP_{y,w}^*)=P_{y,w}^*$ for all $y<w$. Now apply $\hat{\sigma}$ to (2.6). 
This yields the equation
\[T_sC_w= C_{sw}-v^{L(s)}C_w+ \sum_{\atop{y<w}{sy<y}} {\sigma}
(\bM_{y,w}^s) \, C_z \qquad \mbox{if $sw>w$}.\]
Thus, we have $M_{y,w}^s={\sigma}(\bM_{y,w}^s)$ if $sy<y<w<sw$.
Finally, we claim that 
\begin{equation*}
\bM_{y,w}^s \neq 0\quad\Longrightarrow\quad M_{y,w}^s=
\sigma(\bM_{y,w}^s)\neq 0. \tag{$\dagger$}
\end{equation*}
Indeed, if $\bM_{y,w}^s\neq 0$, we write $\bM_{y,w}^s=n_1\gamma_1+\cdots +
n_r \gamma_r$ where $0\neq n_i\in \Z$ and $\gamma_{i-1}^{-1}\gamma_i\in 
\Gamma_+^{(b)}(W)$. By condition ($*$), we have $\sigma(\gamma_{i-1}^{-1}
\gamma_i)=v^{a_i}$ with $a_i>0$ for all~$i$. Consequently, $M_{y,w}^s=\sigma
(\bM_{y,w}^s)$ is a combination of pairwise different powers of~$v$ 
and, hence, non-zero. Thus, ($\dagger$) holds.

So we conclude that two elements satisfy $y\leq_L w$ with respect to 
$\Gamma_+\subset \Gamma$ if and only if they satisfy the analogous 
relation with respect to the weight function $L$. Thus, the relations
$\leq_L$, $\sim_L$, $\leq_{LR}$ and $\sim_{LR}$ are the same in the
two situations, and so are the corresponding representations of~$W$.
\end{proof}

In order to deal with ``distinguished involutions''  as in Theorem~1.2,
we shall need the following remark.

\begin{rem} \label{distinv} In the above setting, let $w\in W$ and write
\begin{equation*} 
\bP_{1,w}^*=\delta_w^{-1}(n_w+\mbox{$\Z$-combination of $\gamma\in 
\Gamma_-$}),\tag{a}
\end{equation*}
where $\delta_w\in \Gamma_+$ and $0\neq n_w\in \Z$. Thus, $\delta_w^{-1}$
is the highest monomial (with respect to the total order specified by
$\Gamma_+\subset \Gamma$) occurring in $\bP_{1,w}^*$. Then $\delta_1=1$
and $\delta_w \in \Gamma_+(W)$ for $w\neq 1$. 

Furthermore, given a left cell $C$ (with respect to $\Gamma_+\subset 
\Gamma$), we write 
\begin{equation*}
\{\delta_w\mid w\in C\}=\{\gamma_1,\gamma_2,\ldots,\gamma_m\} \qquad
\mbox{where $\gamma_{i-1}^{-1}\gamma_i \in \Gamma_+$ for $2\leq i\leq m$}.
\tag{b}
\end{equation*}
Let $\Gamma_+'(W)$ be the union of $\Gamma_+(W)$, the set of all
elements $\gamma^{-1}$ where $\gamma$ occurs in a ${\Z}$-combination 
as in (a) (for any $w\in W$), and the set of all elements
$\gamma_{i-1}^{-1} \gamma_i$ ($2\leq i \leq m$) as in (b) (for any 
left cell $C$ where $m\geq 2$). Assume that 
\begin{equation*}
\sigma(\Gamma_+'(W))\subseteq \{v^n \mid n>0\}.\tag{$*^\prime$}
\end{equation*}
Then, writing $\sigma(\delta_w)=v^{\Delta(w)}$ where $\Delta(w)\in
{\Z}_{\geq 0}$, we have 
\[ P_{1,w}^*=n_wv^{-\Delta(w)}+\mbox{strictly smaller powers of $v$}.\]
Furthermore, if the function $w\mapsto \delta_w$ reaches its mimimum 
at exactly one element in a left cell $C$, then so does the function
$w\mapsto \Delta(w)$.
\end{rem}

\begin{exmp} \label{exp1} Let $W=\langle s,t\rangle$ be a dihedral group 
of order $m\geq 4$, where $m$ is even. Let $v_s$ and $v_t$ be two 
independent indeterminates and consider the ring of Laurent polynomials
$\bA={\Z}[v_s^{\pm 1}, v_t^{\pm 1}]$. Let $\Gamma=\{v_s^iv_t^j
\mid i,j\in \Z\}$ and consider the total order specified by 
\[ \Gamma_+=\{v_s^iv_t^j\mid i>0\} \cup \{v_t^j\mid j>0\}\]
(a lexicographic order where $v_s>v_t$). The polynomials $\bP_{y,w}$
have been determined independently in \cite[Exc.~11.4]{ourbuch} and in
\cite[Chap.~7]{Lusztig03}. Let $y<w$ and write $v_wv_y^{-1}=v_s^{m_s}
v_t^{m_t}$ where $m_s,m_t\geq 0$. Then 
\[ \bP_{y,w}=\left\{\begin{array}{cl} \displaystyle 
\sum_{i=0}^{m_t} (-1)^i v_t^{2i} &\qquad \mbox{if $w<tw$, $w<wt$ 
and $y\leq tsw<sw$}\\
1+v_t^2 & \qquad \mbox{if $w<sw$, $w<ws$ and $y\leq stw<tw$},
\\ 1 & \qquad \mbox{otherwise}.\end{array}\right.\]
The $M$-polynomials are given by
\[ \bM_{y,w}^s=\left\{ \begin{array}{cl} v_sv_t^{-1}+v_s^{-1}v_t & \mbox{if
$l(w)=l(y)+1$, $sy<y<w<sw$}, \\ 1 & \mbox{if $l(w)=l(y)+3$, $sy<y<w<sw$}.
\end{array}\right.\]
All other $M$-polynomials are~$0$. Now consider a weight function $L$ on $W$
such that 
\[ L(s)>L(t)>0.\]
Let $v$ be another indeterminate; then we have a ring homomorphism 
\[\sigma \colon {\Z}[\Gamma]\rightarrow {\Z}[v,v^{-1}],\qquad 
v_s^iv_t^j \mapsto v^{L(s)i+L(t)j}.\] 
We claim that condition ($*$) in Proposition~\ref{prop1} is satisfied. For 
this purpose, we first have to determine the monomials which can occur in
a polynomial $\bP_{y,w}^*$ for $y<w$. Write $v_wv_y^{-1}=v_s^{m_s}v_t^{m_t}$
as above. Since $y<w$, we have $m_s>0$ or $m_t>0$. If $w<tw$, $w<wt$ and 
$y\leq tsw<sw$, then $w$ has a reduced expression which starts and ends 
with $s$. Since $y$ is a subexpression of $w$, we conclude that $m_s
\geq m_t$. Hence $\bP_{y,w}^*$ is a linear combination of monomials 
$v_s^{-m_s}v_t^j$ where $j\leq m_t\leq m_s$. On the other hand, if $w<sw$,
$w<ws$ and $y\leq stw<tw$, then $m_s\geq 1$ and $m_t\geq 1$. So $\bP_{y,w}^*$
is a linear combination of monomials $v_s^{-m_s}v_t^j$ where $j\leq 1$. 
Finally, in the cases where $\bP_{y,w}=1$, we have $\bP_{y,w}^*=
v_s^{-m_s}v_t^{-m_t}$. Thus, we find that 
\[ \Gamma_+(W)\subseteq \{v_s^iv_t^j \mid i\geq 0, i+j\geq 0, (i,j)\neq 
(0,0)\}.\]
Now, if $i\geq 0$ and $i+j\geq 0$, then  $L(s)i+L(t)j\geq L(t)i+L(t)j=
L(t)(i+j)\geq 0$. Furthermore, if $i>0$, then the first inequality is strict
and so $L(s)i+L(t)j>0$; while if $i=0$, then $j>0$ and so $L(s)i+L(t)j>0$. 
Next, we also see that the required condition holds for the monomials
occurring in the polynomials $\bM_{y,w}^t$. Thus, ($*$) holds. 

We conclude that $P_{y,w}^*={\sigma}(\bP_{y,w}^*)$ for all $y<w$ in $W$.
Thus, for any weight function such that $L(s)>L(t)>0$, the corresponding 
polynomials $P_{y,w}^*$ are obtained by specialisation from the polynomials
$\bP_{y,w}^*$ which have been determined for one fixed choice of $\Gamma_+
\subset \Gamma$. Furthermore, the partition of $W$ into left cells is the 
same for all weight functions such that $L(s)>L(t)>0$ (and it is given by 
the partition into left cells with respect to $\Gamma_+\subset \Gamma$). 
An explicit description of these left cells is given in 
\cite[Chap.~8]{Lusztig03}. The distinguished involutions are $1$, $s$, 
$t$, $tst$, $tw_0$, $w_0$. For the left cell representations 
and constructibe representations, see also \cite[\S 6]{myert02}. 
\end{exmp}

\begin{defn} \label{def1} Let $L,L'$ be two weight functions on $W$. We 
say that $L,L'$ are $\Gamma_+$-equivalent if there exists an abelian group 
$\Gamma$, a total order specified by $\Gamma_+\subset \Gamma$ and a set of 
parameters $\{v_s\mid s\in S\}\subset \Gamma_+$ such that the following holds.
\begin{itemize}
\item[(a)] There exist ring homomorphisms $\sigma,\sigma'\colon 
{\Z}[\Gamma]\rightarrow {\Z}[v,v^{-1}]$ such that $\sigma(v_s)=v^{L(s)}$
and $\sigma'(v_s)=v^{L'(s)}$ for all $s\in S$.
\item[(b)] Condition ($*$) in Proposition~\ref{prop1} is satisfied for
both $\sigma$ and $\sigma'$.
\end{itemize}
We say that $L,L'$ are equivalent if $L=L'$ or if there exists a sequence
of weight functions $L=L_0,L_1,\ldots,L_n=L'$ and abelian groups $\Gamma_1,
\ldots,\Gamma_n$ such that $L_{i-1},L_i$ are $(\Gamma_i)_+$-equivalent
with respect to a total order specified by $(\Gamma_i)_+\subset \Gamma_i$
for $1\leq i \leq n$. In particular, equivalent weight functions give rise
to the same partition of $W$ into left cells and to the same set of left 
cell representations of $W$.
\end{defn}

\begin{prop} \label{thm1} Assume that $W$ is finite and let $w_0\in W$ be
the longest element. Then there exists a constant $N \leq 8l(w_0)^3$ such 
that any weight function on $W$ is equivalent to a weight function $L$ such 
that $1 \leq L(s) \leq N$ for all $s\in S$.
\end{prop}

The proof will be given in Section~\ref{sec:proof} (see 
Corollary~\ref{plem3}). Note that, since $W$ is finite, there clearly 
exists some constant $N$ having the above property.  The point about 
Proposition~\ref{thm1} is that we can give an explicit bound for $N$. 
We have not tried to obtain an optimal bound theoretically. However, 
the proofs of Proposition~\ref{plem2} and Corollary~\ref{plem3} will 
show how to determine such a bound efficiently.

\begin{rem} \label{rem1} Let $L\colon W\rightarrow \Z$ be a weight function
such that $L(s)>0$ for all $s\in S$. Let $d>0$ be a positive integer. Then 
the function $L_d\colon W\rightarrow \Z$ defined by $L_d(w):=dL(w)$ also is
a weight function, and we leave it as an (easy) exercise to the reader
to check that $L,L_d$ are equivalent.  Thus, in order to classify 
weight functions up to equivalence, it will be sufficient to consider
only those weight functions $L$ such that $\gcd(\{L(s)\mid s\in S\})=1$.
\end{rem}

\begin{exmp} \label{exp3} In practice, the bound $N$ in Proposition~\ref{thm1}
may be much smaller than $8l(w_0)^3$. For example, if $W=\langle s,t\rangle$ 
is a dihedral group of type $I_2(m)$ (with $m\geq 4$ even), then we may take 
$N=2$. Indeed, let us specify a weight function $L\colon W\rightarrow \Z$ by 
the pair $(a,b)$ such that $L(s)=a$ and $L(t)=b$. Then, by 
Example~\ref{exp1}, there are exactly three equivalence classes 
of weight functions:
\begin{alignat*}{2}
\cL_1 &= \{ (a,b) \mid a=b>0\}, & \qquad 
\mbox{representative:} & \quad (1,1),\\
\cL_2 &= \{ (a,b) \mid a>b>0\}, & \qquad 
\mbox{representative:} & \quad (2,1),\\
\cL_3 &= \{ (a,b) \mid b>a>0\}, & \qquad 
\mbox{representative:} & \quad (1,2).
\end{alignat*}

If $W$ is of type $F_4$, we will see in Section~4 that we may take $N=7$.

Now let $W$ be of type $B_n$, with diagram given as follows.
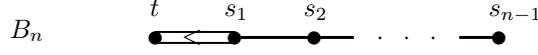
\begin{center}
\begin{picture}(250,30)
\put( 10,10){$B_n$} 
\put( 65,10){\circle*{5}}
\put( 65,12){\line(1,0){30}}
\put( 65, 8){\line(1,0){30}}
\put( 75,7.5){$<$}
\put( 95,10){\circle*{5}}
\put( 95,10){\line(1,0){30}}
\put(125,10){\circle*{5}}
\put(125,10){\line(1,0){15}}
\put(150,10){\circle*{1}}
\put(160,10){\circle*{1}}
\put(170,10){\circle*{1}}
\put(180,10){\line(1,0){15}}
\put(195,10){\circle*{5}}
\put( 63,18){$t$}
\put( 91,18){$s_1$}
\put(121,18){$s_2$}
\put(191,18){$s_{n-1}$}
\end{picture}
\end{center}
Here, the generators $s_i$ are all conjugate, while $t$ and $s_1$ are not
conjugate. Thus, a weight function $L\colon W\rightarrow \Z$ is 
uniquely specified by the values
\[ b:=L(t)>0 \qquad \mbox{and}\qquad a:=L(s_1)=L(s_2)=\cdots =L(s_{n-1})>0.\]
The best bound does not yet seem to be known. Recently, Bonnaf\'e and Iancu
have shown that all weight functions such that $a/b>n-1$ are equivalent.
Experiments with {\sf CHEVIE} lead to the following general conjecture.
\end{exmp}

\begin{conj} \label{conj1} In type $B_n$ with diagram and weight 
function as specified above, we have the following equivalence
classes of weight functions.
\begin{align*}
\cL_0 &= \{ (b,a,a,\ldots,a) \mid a>b>0\}, \\
\cL_1 &= \{ (a,a,a,\ldots,a) \mid a>0\} \quad 
(\mbox{equal parameter case}),\\
\cL_i &= \{ (ia,a,a,\ldots,a) \mid a>0\} 
\;\;(\mbox{where $2\leq i \leq n-1$}), \\
\cL_{i,i-1} & =\{(b,a,a,\ldots,a) \mid ia>b>(i-1)a>0\}
\;\;(\mbox{where $2\leq i \leq n-1$}), \\
\cL_{\text{asymp}} &= \{ (b,a,a,\ldots,a) \mid b>(n-1)a>0\}.
\end{align*}
(The functions in $\cL_{\text{asymp}}$ 
correspond to the case treated by Bonnaf\'e--Iancu \cite{BI}.)
\end{conj}

Furthermore, if \S 2 (C) holds, then all left cell representations with
respect to $L$ will be irreducible, unless we have $L\in \cL_i$ for
some $1\leq i \leq n-1$ (see \cite[22.25]{Lusztig03}); if $L\in \cL_i$
for some $i$, then the left cell representations will be given as in 
\cite[22.24]{Lusztig03}. 

The above conjecture is a slightly different version (via the exact
form of the equivalence relation in Definition~\ref{def1}) of a part 
of several conjectures that were formulated by Bonnaf\'e (private
communication).  We will see in Section~4 (where we consider $W$ of 
type $F_4$) that weight functions which are not equivalent may still 
give rise to the same partition of $W$ into left cells. This phenomenon 
does not seem to occur in type $B_n$.

Using our {\sf CHEVIE}-program, we have verified the above conjecture 
for $B_3$ and $B_4$. For example, in type $B_4$, we obtain $7$ equivalence 
classes of weight functions:
\begin{alignat*}{2}
\cL_0 &= \{ (b,a,a,a) \mid a>b>0\}, 
& \quad \mbox{representative:} & \quad (1,2,2,2),\\
\cL_1 &= \{ (a,a,a,a) \mid a>0\}, 
& \quad \mbox{representative:} & \quad (1,1,1,1),\\
\cL_2 &= \{ (2a,a,a,a) \mid a>0\}, 
& \quad \mbox{representative:} & \quad (2,1,1,1),\\
\cL_3 &= \{ (3a,a,a,a) \mid a>0\}, 
& \quad \mbox{representative:} & \quad (3,1,1,1),\\
\cL_{2,1} &= \{ (b,a,a,a) \mid 2a>b>a>0\}, 
& \quad \mbox{representative:} & \quad (3,2,2,2),\\
\cL_{3,2} &= \{ (b,a,a,a) \mid 3a>b>2a>0\}, 
& \quad \mbox{representative:} & \quad (5,2,2,2),\\
\cL_{\text{asymp}} &= \{ (b,a,a,a) \mid b>3a>0\}, 
& \quad \mbox{representative:} & \quad (4,1,1,1).
\end{alignat*}



The above results are only concerned with finite Coxeter groups. It would 
be interesting to study equivalence classes of weight functions for affine
Weyl groups.

\section{On the equivalence classes of weight functions} \label{sec:proof} 

We place ourselves in the general setting where $W$ is any Coxeter group
with generators $S$ and where we are given an abelian group $\Gamma$, a 
total order specified by $\Gamma_+\subset \Gamma$ and a set of parameters 
$\{v_s\mid s\in S\}\subset \Gamma_+$ for the corresponding 
Iwahori--Hecke algebra of $W$. One of the aims of this section is to
provide a proof of Proposition~\ref{thm1}. Our first task will be to get 
some control on the degrees of the monomials that might occur in the 
polynomials $\bP_{y,w}^*$ and $\bM_{y,w}^s$. Now, Lusztig gives some rather 
explicit such bounds, but only in the setting involving a weight function, 
and these are not entirely sufficient for our purposes. To illustrate our 
point, consider the following example.

\begin{exmp} \label{exp2}  Let $W=\langle s,t\rangle$ be a dihedral
group as in Example~\ref{exp1}. Consider a weight function $L$ where
$L(s)=a>1$ is a big number and $L(t)=1$. Then \cite[Prop.~6.4]{Lusztig03} 
tells us that $M_{y,w}^s$ is a $\Z$-linear combination of powers $v^n$ with
$-a+1\leq n\leq a-1$ and $n\equiv L(w)-L(y)-L(s) \bmod 2$. So, a priori, 
$M_{y,w}^s$ could be a polynomial involving many non-zero terms. However, 
from the formula given in Example~\ref{exp1} and Proposition~\ref{prop1}, 
we see that $M_{y,w}^s$ only involves very few terms:
\[ M_{y,w}^s={\sigma}(\bM_{y,w}^s)=\left\{ \begin{array}{cl} 
v^{a-1}+v^{1-a} & \mbox{if $l(w)=l(y)+1$, $sy<y<w<sw$}, \\ 1 & 
\mbox{if $l(w)=l(y)+3$, $sy<y<w<sw$}.  \end{array}\right.\]
To explain this behaviour, we need to establish some bounds in the general
framework with respect to an abelian group $\Gamma$ and a total order on it.
\end{exmp}

\begin{lem} \label{lem3} Let $y,w \in W$ be such that $y \leq w$. Then
the following hold.
\begin{itemize}
\item[(a)] $v_wv_y^{-1}\bP_{y,w}^*$ is a polynomial in $\{v_s^2\mid 
s \in S\}$, with constant term~$1$.
\item[(b)] $v_wv_y^{-1}\overline{\bP}_{y,w}^*$ is a polynomial in 
$\{v_s^2\mid s\in S\}$, with constant term~$0$.
\end{itemize}
\end{lem}

\begin{proof} The following proof is more or less a copy of that of 
\cite[Prop.~5.4]{Lusztig03}. However, in \cite{Lusztig03}, Lusztig 
exclusively considers the situation involving a weight function. Thus,
in order to show that all the arguements go through in the general case,
we include the details here. First, we shall need the $R$-polynomials in
the general setting, as defined in \cite{Lusztig83}. For $y\in W$, we have 
\[ \overline{\bT}_y=\bT_{y^{-1}}^{-1}=\sum_{x \in W}\overline{\bR}_{x,y} 
\bT_x \qquad \mbox{where $\bR_{x,y}\in {\Z}[\Gamma]$}.\]
We have the following recursion formula. If $sy<y$ for some $s \in S$, then
\begin{alignat*}{2}
\bR_{x,y}&=\bR_{sx,sy}+(v_s-v_s^{-1})\bR_{x,sy}&&\qquad\mbox{if $sx>x$},\\
\bR_{x,y}&=\bR_{sx,sy} &&\qquad \mbox{if $sx<x$}.
\end{alignat*}
(Same proof as in \cite[Lemma~4.4]{Lusztig03}.) Using the above recursion 
formula, one easily shows that $\bR_{y,y}=1$ and $\bR_{x,y}=0$ unless 
$x \leq y$. Furthermore, 
\begin{equation*}
v_yv_x^{-1}\bR_{x,y} \in {\Z}[v_s^2\mid s\in S], \mbox{ with constant 
term~$(-1)^{l(y)-l(x)}$}.\tag{$*$}
\end{equation*}
(Same proof as in \cite[Lemma~4.7]{Lusztig03}.) The Kazhdan--Lusztig 
polynomials and the $R$-polynomials are related by the following identity 
(see \cite[Prop.~2]{Lusztig83}). We have 
\[\overline{\bP}^{*}_{x,w}-\bP^{*}_{x,w} = \sum_{x<y\leq w} 
\bR_{x,y}\bP^{*}_{y,w} \qquad \mbox{for all $x<w$ in $W$}.\]
Now, for the proof of (a) and (b), we proceed by induction on $l(w)-l(y)$. 
If $y=w$, then $\bP_{w,w}^*=1$ and there is nothing to prove. Now assume 
that $y<w$. Multiplying both sides of the identity relating 
Kazhdan--Lusztig polynomials and $R$-polynomials with $v_w v_y^{-1}$ yields
\[ v_wv_y^{-1}\overline{\bP}^{*}_{y,w}-v_wv_y^{-1}\bP^{*}_{y,w}=
\sum_{y <x \leq w} (v_xv_y^{-1}\bR_{y,x})(v_wv_x^{-1}\bP^{*}_{x,w}).\]
By induction and ($*$), all terms on the right hand side are polynomials 
in $\{v_s^2\mid s \in S\}$. Hence so is the left hand side. Since 
$\bP^*_{y,w}$ and $\overline{\bP}^*_{y,w}$ have no terms in common, we 
conclude that both $v_wv_y^{-1}\bP_{y,w}^*$ and $v_wv_y^{-1} 
\overline{\bP}_{y,w}^*$ are polynomials in the variables $v_s^2$ 
($s\in S$). Now consider the constant terms on both sides of the above 
equation. We begin with the right hand side. By induction and ($*$), it 
has constant term
\[ \sum_{y<x\leq w}(-1)^{l(x)-l(y)} \cdot 1=-1+ (-1)^{l(y)}\sum_{y
\leq x\leq w}(-1)^{l(x)} =-1,\] 
where the last equality holds by \cite[Prop.~4.8]{Lusztig03} (an identity
due to D.~N.~Verma). It remains to observe that $v_wv_y^{-1}
\overline{\bP}^*_{y,w}\in {\Z}[\Gamma_+]$ and so the constant term is~$0$.
Hence the constant term of $-v_wv_y^{-1}\bP_{y,w}^*$ equals~$-1$, as 
required.
\end{proof}

\begin{lem} \label{lem4} Let $y,w \in W$ and $s \in S$ be such that
$sy<y<w<sw$. Then $v_sv_wv_y^{-1} \bM_{y,w}^s$ is a polynomial in
$\{v_t^2\mid t \in S\}$, with constant term~$0$.
\end{lem}

\begin{proof} As in the proof of \cite[Prop.~4]{Lusztig83}, one
considers the identity (arising from (2.6)):
\[\bT_s\bC_w-\bC_{sw}+v_s\bC_w-\sum_{\atop{y<w}{sy<y}} \bM_{y,w}^s \, 
\bC_z=0.\]
Expressing all terms in the basis $\{\bT_y \mid y\in W\}$ of $\bH$, the
coefficient of every $\bT_y$ must be zero. That coefficient is given by
\[ f_y=v_s\bP_{y,w}^*+\bP_{sy,w}^*-\bP_{y,sw}^*-\sum_{\atop{y\leq z<w}{sz<z}}
\bP_{y,z}^*\bM_{z,w}^s.\]
Hence, given that $f_y=0$, we obtain
\[ \bM_{y,w}^s=\bP^{*}_{sy,w} -\bP^{*}_{y,sw}+v_s\bP^{*}_{y,w}-
\sum_{\atop{y<z<w}{sz<z}} \bP^{*}_{y,z}\bM_{z,w}^s.\] 
Since $sy<y$ and $sw>w$, we
have $v_y=v_sv_{sy}$ and $v_{sw}=v_sv_w$. Thus, multiplying the above
equation by $v_sv_wv_y^{-1}$ yields that
\[ v_sv_wv_y^{-1}\bM_{y,w}^s=\bP_{sy,w} -\bP_{y,sw}+v_s^2\bP_{y,w}-
\sum_{\atop{y <z<w}{sz<z}} \bP_{y,z}(v_sz_wv_z^{-1}\bM_{z,w}^s).\]
Hence, the assertion follows by induction on $l(w)-l(y)$ and using 
Lemma~\ref{lem3}.
\end{proof}

{From} now on, we assume that $W$ is finite and let $w_0\in W$ be the 
longest element. Then, by the classification of
finite Coxeter groups, unequal parameters can only occur for $W$ of type
$I_2(m)$ (with $m$ even), $B_n$ (any $n\geq 3$) or $F_4$. Furthermore,
in these cases, a weight function on $W$ may take at most $2$ different
values on the generators of $W$. Thus, we will now consider an abelian 
group $\Gamma=\{x^iy^j \mid i,j\in \Z\}$ where $x$ and $y$ are independent 
indeterminates and where $\Gamma_+\subset \Gamma$ is any total order.
Furthermore, let $S=S_x \amalg S_y$ be a partition (where
$S_x,S_y\neq \varnothing$) such that no generator in $S_x$ is conjugate 
to any generator in $S_y$. The parameters of the corresponding 
Iwahori--Hecke algebra will be assumed to be given by 
\[ v_s=x \quad (\mbox{if $s\in S_x$}) \qquad \mbox{and}\qquad v_t=y
\qquad (\mbox{if $t\in S_y$}).\]

\begin{lem} \label{lem5} The monomials involved in any polynomial
$\bP_{y,w}^*$ or in any polynomial $\bM_{y,w}^s$ are of the form
$x^iy^j$ where $-l(w_0)<i,j<l(w_0)$. In particular, we have
$\Gamma_+(W)\subseteq \{x^iy^j\mid  -l(w_0)<i,j<l(w_0)\}$.
\end{lem}

\begin{proof} Let $y,w\in W$, $y\leq w$. Thus, since $y$ is a subexpression
of $w$, we have $v_wv_y^{-1}=x^ay^b$ where $a,b\geq 0$. Furthermore, let us 
write $\bP_{y,w}^*=\sum_{(i,j) \in I} n_{ij} x^iy^j$ where $I \subseteq 
\Z \times \Z$ is a finite subset and $n_{ij} \in \Z$. Thus, using 
Lemma~\ref{lem3}, we have 
\begin{align*}
v_wv_y^{-1}\bP_{y,w}^* &=\sum_{(i,j)\in I} n_{ij} x^{a+i}y^{b+j} \in 
{\Z}[x^2,y^2],\\
v_wv_y^{-1}\overline{\bP}^*_{y,w} &=\sum_{(i,j)\in I} n_{ij} 
x^{a-i}y^{b-j} \in {\Z}[x^2,y^2].
\end{align*}
Now let $(i,j)\in I$. We certainly have $0\leq a,b<l(w_0)$. This
yields $0 \leq a+i<l(w_0)+i$ and $0\leq a-i <l(w_0)-i$. 
Consequently, we have $-l(w_0)<i<l(w_0)$.  A similar argument 
shows that we also have $-l(w_0)< j < l(w_0)$.

Now assume that $sy<y<w<sw$ and write $\bM_{y,w}^s=f+c+\overline{f}$
where $c\in \Z$ and $f\in {\Z}[x^{\pm 1},y^{\pm 1}]$. Let 
$f=\sum_{(i,j)\in J} f_{ij} x^iy^j$ where $J\subseteq \Z\times \Z$ is
a finite subset and $f_{ij}\in \Z$. As above, we see that 
$v_sv_wv_y^{-1}= x^ay^b$ where $0\leq a,b<l(w_0)$. (Note that
$y<w<w_0$.) Using Lemma~\ref{lem5}, this yields
\begin{align*}
v_sv_wv_y^{-1}&\bM_{y,w}^s=x^ay^bf+cx^ay^b+x^ay^b\overline{f}\\
&= cx^ay^b+\sum_{(i,j)\in J}f_{ij}(x^{a+i}y^{b+j}+x^{a-i}y^{b-j})
\in {\Z}[x^2,y^2].
\end{align*}
Arguing as above, we see that $-l(w_0)<i,j<l(w_0)$ for all
$(i,j)\in J$.
\end{proof}

Now, a weight function $L\colon W\rightarrow \Z$ is uniquely specified by
the values
\[ a:=L(s)>0 \quad \mbox{(where $s\in S_x$)} \quad \mbox{and} \quad
b:=L(t)>0 \quad \mbox{(where $t\in S_y$)}.\]
We shall write $L=L_{a,b}$. Let us consider the set 
\[ \cE:=\{x\in {\Q}_{>0}\mid x=\pm i/j \mbox{ where } i,j \neq 0 
\mbox{ and } -2l(w_0)<i,j<2l(w_0)\}\]
and write $\cE=\{x_1,\ldots,x_n\}$ where $0<x_1<x_2<\cdots <x_n$.  By
convention, we set $x_0=0$ and $x_{n+1}=\infty$. For any $0 \leq k\leq n$, 
we consider the set of weight functions 
\[\cL_k :=\{L_{a,b} \mid a,b>0 \mbox{ such that } x_{k} <b/a <x_{k+1}\}.\]
Let us fix $0\leq k \leq n$ and write $x_k=d/c$ where $c,d$ are integers
such that $0\leq c,d<2l(w_0)$ and $c\neq 0$. Then we consider the total 
order in $\Gamma$ specified by 
\[\Gamma_+^{(k)} = \{x^iy^j \mid ci+dj>0\}\cup \{x^iy^j\mid ci+dj=0 
\mbox{ and } i>0\} \quad \mbox{if $d\geq c$},\]
or 
\[ \Gamma_+^{(k)} = \{x^iy^j \mid ci+dj>0\}\cup \{x^iy^j\mid ci+dj=0 
\mbox{ and } j>0\} \quad \mbox{if $d<c$};\]
(a weighted lexicographic order). Note that, if $k=d=0$, then 
\[ \Gamma_+^{(0)}=\{x^iy^j \mid i>0,j\in \Z\}\cup \{y^j\mid j>0\}; \]
(a pure lexicographic order). 

\begin{prop} \label{plem2} In the above settung, all the weight functions
in $\cL_k$ are $\Gamma_+^{(k)}$-equivalent.
\end{prop}

\begin{proof} Let $a,b>0$ be such that $x_k<b/a< x_{k+1}$. The idea
is to get some control on the set $\Gamma_+(W) \subseteq \Gamma_+$ and
to show that condition ($*$) in Proposition~\ref{prop1} is satisfied for 
the ring homomorphism
\[ \sigma_{a,b} \colon {\Z}[\Gamma] \rightarrow {\Z}[v,v^{-1}],
\qquad x^iy^j \mapsto v^{ai+bj}\]
and the total order $\Gamma_{+}^{(k)}\subset \Gamma$ specified above.
Now, by Lemma~\ref{lem5}, we have
\[ \Gamma_+(W)\subseteq \{x^iy^j \mid x^iy^j \in \Gamma_+ \mbox{ and }
-2l(w_0)<i,j<2l(w_0)\}.\]
To check condition ($*$), assume first that $c<d$. Let $x^iy^j\in 
\Gamma_+^{(k)}(W)$. In particular, this means that $ci+dj\geq 0$. 
Furthermore, we have $-2l(w_0)<i,j<2l(w_0)$ and so $\pm i/j \in \cE$. 
Now, we must show that $ai+bj>0$. If $i=0$ or $j=0$, this is clear. If 
$j>0$, then we have 
\[ ai+bj=a(i+jb/a)>a(i+x_kj)=a(i+jd/c)=(a/c)(ci+dj) \geq 0,\]
as required. Next assume that $j<0$. Then, by the definition of 
$\Gamma_+^{(k)}$ (recall that we are assuming $c<d$), we must have 
$ci+dj>0$ and so $-i/j>d/c=x_k$. Now, if we had $ai+bj\leq 0$, then 
we would obtain 
\[ x_k <-i/j \leq b/a<x_{k+1}\]
and so $-i/j\not\in \cE$, a contradiction. Thus, condition ($*$) holds.
The argument for the case where $d\leq c$ is completely analogous.
\end{proof}


\begin{cor} \label{plem3} Let $\cE=\{x_1,\ldots,x_n\}$ as above.
Let $L=L_{a,b}$ be any weight function on $W$ where $a,b>0$.
\begin{itemize}
\item[(1)] If $b/a=x_k$ for some $1\leq k\leq n$, then $L_{a,b}$ is 
equivalent to $L_{c,d}$ where $0<c,d<2l(w_0)$ are such that $b/a=d/c$.
\item[(2)] If $b/a\not\in \cE$, then there exist integers $1\leq a',b'
\leq 8l(w_0)^3$ such that $L_{a,b}$ is equivalent to $L_{a',b'}$.
\end{itemize}
\end{cor}

\begin{proof} Recall that $x_0=0$ and $x_{n+1}=\infty$. Hence
there exists some $k\in \{0,1,\ldots,n\}$ such that $x_{k} \leq b/a
<x_{k+1}$. We write $x_{k}=d/c$ where $0\leq c,d<2l(w_0)$ and 
$c\neq 0$.  If $x_{k}=b/a$, then $L_{a,b},L_{c,d}$ are equivalent by 
Remark~\ref{rem1}. Thus, (1) is proved. Now assume that $x_{k}<b/a<
x_{k+1}$. Since both $x_k$ and $x_{k+1}$ are rational numbers where the 
numerator and the denominator are strictly bounded by $2l(w_0)$, we 
certainly have $1/4l(w_0)^2<x_{k+1}-x_{k}$.
Furthermore, note that $x_n<2l(w_0)$. Thus, we can find some integers 
$a',b'$ such that $1\leq a',b'\leq 8l(w_0)^3$ and $x_k<b'/a'<x_{k+1}$. 
Then $L_{a,b}$ and $L_{a',b'}$ are equivalent by Proposition~\ref{plem2}.
Thus, (2) is proved.
\end{proof}

\begin{exmp} \label{expasy} Assume that $a,b>0$ are such that $a/b 
\geq 2l(w_0)$. Then $L_{a,b}$ is equivalent to the weight function 
$L_{2l(w_0),1}$.

To see this, note that $1/2l(w_0)<x_1$. Hence, we are in the case where 
$b/a \leq 1/2l(w_0)<x_1$. Thus, we have $L_{a,b} \in \cL_0$. By 
Proposition~\ref{plem2}, all weight functions in $\cL_0$ are equivalent. 
It remains to note that $L_{2l(w_0),1}$ also belongs to $\cL_0$. 

This example provides a more formal justification for \cite[Remark 6.1]{BI}.
\end{exmp}

\section{Kazhdan--Lusztig polynomials and left cells in type $F_4$} 
\label{sec:f4}
Our aim is to work out the equivalence classes of weight functions on a 
Coxeter group of type $F_4$. Throughout this section, let $W$ be
a Coxeter group of type $F_4$, with generating set $S=\{s_1,s_2,s_3,s_4\}$
and Dynkin diagram given as follows.
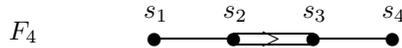
\begin{center}
\begin{picture}(200,30)
\put( 10,10){$F_4$} 
\put( 61,18){$s_1$}
\put( 91,18){$s_2$}
\put(121,18){$s_3$}
\put(151,18){$s_4$}
\put( 65,10){\circle*{5}}
\put( 95,10){\circle*{5}}
\put(125,10){\circle*{5}}
\put(155,10){\circle*{5}}
\put(105,7.5){$>$}
\put( 65,10){\line(1,0){30}}
\put( 95,12){\line(1,0){30}}
\put( 95, 8){\line(1,0){30}}
\put(125,10){\line(1,0){30}}
\end{picture}
\end{center}
There are $25$ irreducible representations of $W$, denoted by 
\[ 1_1, 1_2,1_3,1_4, 2_1, 2_2, 2_3, 2_4, 4_1, 4_2, 4_3, 4_4, 4_5, 6_1,
6_2, 8_1, 8_2, 8_3, 8_4, 9_1, 9_2,  9_3, 9_4, 12_1, 16_1;\] 
see \cite[4.10]{LuBook} or \cite[5.3.6 and Table~C.3]{ourbuch}. The 
generators $s_1,s_2$ are conjugate in $W$, and so are the generators 
$s_3,s_4$ (while $s_2$ and $s_3$ are not conjugate). Thus, a weight 
function $L\colon W\rightarrow \Z$ is uniquely determined by
\[ L(s_1)=L(s_2)=a>0\qquad \mbox{and} \qquad L(s_3)=L(s_4)=b>0.\]
We shall denote such a weight function by $L=L_{a,b}$. By the symmetry
of the above diagram, we may assume throughout that $a\leq b$. 

Let $x,y$ be independent indeterminates over $\Z$ and consider the abelian 
group
\[ \Gamma=\{x^iy^j\mid i,j\in \Z\}.\]
Let $v$ be another indeterminate. Then we have a ring homomorphism
\[ \sigma_{a,b} \colon {\Z}[\Gamma] \rightarrow {\Z}[v,v^{-1}],\qquad
x^iy^j \mapsto v^{ai+bj}.\]
Now, in type $F_4$, we have $l(w_0)=24$ and so, by Corollary~\ref{plem3},
we know that $L_{a,b}$ is equivalent to a weight function $L_{c,d}$ where 
$1\leq c\leq d\leq 48^3=110592$. In principle, we could just go through 
all these possibilities, determine the corresponding left cell 
representations and so on---but these are far too many cases! However, 
now we can use our {\sf CHEVIE}-program to compute explicitly all the 
polynomials $\bP_{y,w}^*$ and $\bM_{y,w}^s$ for any total order on $\Gamma$. 
The explicit knowledge of these polynomials will yield much sharper bounds 
than the general bounds obtained in Lemma~\ref{lem5}. 

As a first illustration of this idea, we consider the following case.

\begin{lem} \label{lem41} Consider the total order on $\Gamma$ specified by 
\[ \Gamma_+=\{x^iy^j \mid j>0,i\in \Z\}\cup \{x^i\mid i>0\}.\]
Then condition ($*$) in Proposition~\ref{prop1} is satisfied for all 
weight functions $L_{a,b}$ such that $b/a>4$. In particular, all
these weight functions are $\Gamma_+$-equivalent. 
\end{lem}

\begin{proof} The idea is basically the same as in the proof of 
Proposition~\ref{plem2}. In fact, the general strategy in 
Corollary~\ref{plem3} shows that all $L_{a,b}$ are $\Gamma_{+}$-equivalent, 
provided that $b/a>2l(w_0)=48$. But now we use our {\sf CHEVIE}-program 
to compute explicitly all the polynomials $\bP_{y,w}^*$ and $\bM_{y,w}^*$ 
(with respect to $\Gamma_+\subset \Gamma$). By inspection of all these
polynomials, we find that
\[\Gamma_+(W)\subseteq\{x^i\mid i>0\}\cup\{x^iy^j\mid j>0,i+4j\geq 0\}.\]
Now let us check that condition ($*$) in Proposition~\ref{prop1} holds
for $\sigma_{a,b}$ provided that $b>4a$. Let $i,j\in \Z$ be such that 
$x^iy^j\in \Gamma_+(W)$. We must show that $ai+bj>0$. If $j=0$, then 
$i>0$ and so $ai+bj=ai>0$. On the other hand, if $j>0$ and $i+4j\geq 0$, 
then $ai+bj=a(i+jb/a)>a(i+4j)\geq 0$, as required. 

We can now apply 
Proposition~\ref{prop1} and conclude that all weight functions $L_{a,b}$
such that $b/a>4$ are $\Gamma_+$-equivalent.
\end{proof}

In order to deal with weight functions $L_{a,b}$ such that $b/a<4$, we
now proceed as follows. We look again at the elements in $\Gamma_+(W)$
computed in the proof of Lemma~\ref{lem41}. Let
\[ \cE=\{x\in {\Q}_{>0} \mid x=\pm i/j \mbox{ where } j\neq 0,x^iy^j\in 
\Gamma_+(W)\}.\]
Then we note that the largest element of $\cE$ below $4$ is $3$. This
leads us to consider weight functions $L_{a',b'}$ where $b'/a'>3$.

\begin{lem} \label{lem42} Consider the total order on $\Gamma$ specified by 
\[ \Gamma_+=\{x^iy^j \mid i+3j>0\}\cup \{x^{-3j}y^j\mid j>0\}.\]
Then condition ($*$) in Proposition~\ref{prop1} is satisfied for all 
weight functions $L_{a,b}$ such that $4>b/a>3$. In particular, all these
weight functions are $\Gamma_+$-equivalent. 
\end{lem}

\begin{proof} This is completely analogous to that of Lemma~\ref{lem41}.
Now we find that 
\begin{multline*}
\Gamma_+(W)\subseteq \{x^i\mid i>0\} \cup\{x^iy^j\mid j>0,i+j\geq 0\}\\
\cup \{x^iy^j\mid i>-j>0,-i/j\geq 4\}\cup \{x^iy^j\mid -i>j>0,-i/j\leq 3\}.
\end{multline*} 
As before, we see that condition ($*$) in Proposition~\ref{prop1} holds, 
provided that $4a>b>3a$. Indeed, let $i,j$ be such that $x^iy^j\in
\Gamma_+(W)$. If $j=0$, then $i>0$ and so $ai+bj=ai>0$. If $j>0$ and
$i+j\geq 0$, then $ai+bj>ai+3aj>a(i+j)\geq 0$. If $i>-j>0$ and $-i/j\geq 4$, 
then $ai+bj=i(a+bj/i)>ia(1+4j/i)\geq 0$. Finally, if $-i>j>0$ and $-i/j>3$, 
then $ai+bj=j(ai/j+b)>aj(i/j+3)\geq 0$, as required. 
\end{proof}

As before, we now look again at the elements in $\Gamma_+(W)$ computed 
in the proof of Lemma~\ref{lem42}. Define $\cE$ in a similar way as above.
Then we note that the largest element of $\cE$ below $3$ is $5/2$. This
leads us to the following case.

\begin{lem} \label{lem43} Consider the total order on $\Gamma$ specified by 
\[ \Gamma_+=\{x^iy^j \mid 2i+5j>0\}\cup \{x^{-5j}y^{2j}\mid j>0\}.\]
Then condition ($*$) in Proposition~\ref{prop1} is satisfied for all 
weight functions $L_{a,b}$ such that $3>b/a>5/2$. In particular, all these
weight functions are $\Gamma_+$-equivalent. 
\end{lem}

\begin{proof} Again, this is completely analogous to that of 
Lemma~\ref{lem41}. Now we find that
\begin{multline*}
\Gamma_+(W)\subseteq \{x^i\mid i>0\} \cup\{x^iy^j\mid j>0,i+j\geq 0\}\\
\cup \{x^iy^j\mid i>0,i+3j\geq 0\} \cup \{x^iy^j\mid -i>j>0,-i/j\leq 5/2\}.
\end{multline*} 
We omit further details.
\end{proof}

We now continue the above procedure. This yields the following cases.

\begin{lem} \label{lem44} Consider the total order on $\Gamma$ specified by 
\[ \Gamma_+=\{x^iy^j \mid i+2j>0\}\cup \{x^{-2j}y^{j}\mid j>0\}.\]
Then condition ($*$) in Proposition~\ref{prop1} is satisfied for all 
weight functions $L_{a,b}$ such that $5/2>b/a>2$. In particular, all these
weight functions are $\Gamma_+$-equivalent. 
\end{lem}

\begin{proof} Again, this is completely analogous to that of 
Lemma~\ref{lem41}. Now we find that 
\begin{multline*}
\Gamma_+(W)\subseteq \{x^i\mid i>0\}\cup\{x^iy^j\mid j>0,i+j\geq 0\}\\
\cup \{x^iy^j\mid i>-j>0,-i/j\geq 5/2\}\cup \{x^iy^j\mid -i>j>0,-i/j\leq 2\}.
\end{multline*} 
We omit further details.
\end{proof}

%

\begin{lem} \label{lem45} Consider the total order on $\Gamma$ specified by 
\[ \Gamma_+=\{x^iy^j \mid 2i+3j>0\}\cup \{x^{-3j}y^{2j}\mid j>0\}.\]
Then condition ($*$) in Proposition~\ref{prop1} is satisfied for all 
weight functions $L_{a,b}$ such that $2>b/a>3/2$. In particular, all these
weight functions are $\Gamma_+$-equivalent. 
\end{lem}

\begin{proof} Again, this is completely analogous to that of 
Lemma~\ref{lem41}. Now we find that 
\begin{multline*}
\Gamma_+(W)\subseteq \{x^i\mid i>0\} \cup\{x^iy^j\mid j>0,i+j\geq 0\}\\
\cup \{x^iy^j\mid i>-j>0,-i/j\geq 2\} \cup 
\{x^iy^j\mid -i>j>0,-i/j\leq 3/2\}.
\end{multline*} 
We omit further details.
\end{proof}

\begin{lem} \label{lem46} Consider the total order on $\Gamma$ specified by 
\[ \Gamma_{+}=\{x^iy^j \mid 3i+4j>0\}\cup \{x^{-4j}y^{3j}\mid j>0\}.\]
Then condition ($*$) in Proposition~\ref{prop1} is satisfied for all 
weight functions $L_{a,b}$ such that $3/2>b/a>4/3$. In particular, all these
weight functions are $\Gamma_+$-equivalent. 
\end{lem}

\begin{proof} Again, this is completely analogous to that of 
Lemma~\ref{lem41}. Now we find that 
\begin{multline*}
\Gamma_+(W)\subseteq \{x^i\mid i>0\} \cup\{x^iy^j\mid j>0,i+j\geq 0\}\\
\cup \{x^iy^j\mid i>-j>0,-i/j\geq 3/2\} \cup \{x^{-4j}y^{3j}\mid j>0\}.
\end{multline*} 
We omit further details.
\end{proof}

\begin{lem} \label{lem47} Consider the total order on $\Gamma$ specified by 
\[ \Gamma_{+}=\{x^iy^j \mid i+j>0\}\cup \{x^{-j}y^{j}\mid j>0\}.\]
Then condition ($*$) in Proposition~\ref{prop1} is satisfied for all 
weight functions $L_{a,b}$ such that $4/3>b/a>1$. In particular, all these
weight functions are $\Gamma_+$-equivalent. 
\end{lem}

\begin{proof} Again, this is completely analogous to that of 
Lemma~\ref{lem41}. Now we find that 
\[\Gamma_+(W)\subseteq \{x^i\mid i>0\} \cup\{x^iy^j\mid j>0,i+j\geq 0\}
\cup \{x^iy^j\mid i>0,3i+4j\geq 0\}.\]
We omit further details.
\end{proof}

\begin{table}[htbp]
\caption{Partial order on two-sided cells in type $F_4$} 
\label{lrgraph} 
\begin{center}
\begin{picture}(325,265)
\put( 16,  0){$a=b$}
\put( 27, 20){\circle{4}}
\put( 32, 18){\scriptsize{$1_4$}}
\put( 27, 22){\line(0,1){16}}
\put( 27, 40){\circle{4}}
\put( 32, 38){\scriptsize{\fbox{$4_5$}}}
\put( 27, 42){\line(0,1){16}}
\put( 27, 60){\circle{4}}
\put( 32, 58){\scriptsize{$9_4$}}
\put( 28, 62){\line(4,5){13}}
\put( 26, 62){\line(-4,5){13}}
\put( 12, 80){\circle{4}}
\put(  0, 78){\scriptsize{$8_2$}}
\put( 42, 80){\circle{4}}
\put( 47, 78){\scriptsize{$8_4$}}
\put( 13, 82){\line(4,5){13}}
\put( 41, 82){\line(-4,5){13}}
\put( 27,100){\circle{4}}
\put(  0, 99){\scriptsize{\fbox{$12_1$}}}
\put( 28,102){\line(4,5){13}}
\put( 26,102){\line(-4,5){13}}
\put( 12,120){\circle{4}}
\put(  0,118){\scriptsize{$8_1$}}
\put( 42,120){\circle{4}}
\put( 47,118){\scriptsize{$8_3$}}
\put( 13,122){\line(4,5){13}}
\put( 41,122){\line(-4,5){13}}
\put( 27,140){\circle{4}}
\put( 32,140){\scriptsize{$9_1$}}
\put( 27,142){\line(0,1){16}}
\put( 27,160){\circle{4}}
\put( 32,158){\scriptsize{\fbox{$4_2$}}}
\put( 27,162){\line(0,1){16}}
\put( 27,180){\circle{4}}
\put( 32,178){\scriptsize{$1_1$}}
\put( 93,  0){$b=2a$}
\put(105, 20){\circle{4}}
\put(110, 18){\scriptsize{$1_4$}}
\put(105, 22){\line(0,1){16}}
\put(105, 40){\circle{4}}
\put(110, 38){\scriptsize{$2_4$}}
\put(105, 42){\line(0,1){16}}
\put(105, 60){\circle{4}}
\put(110, 58){\scriptsize{$4_5$}}
\put(105, 62){\line(0,1){16}}
\put(105, 80){\circle{4}}
\put( 84, 78){\scriptsize{\fbox{$1_2$}}}
\put(106, 82){\line(1,1){14}}
\put(104, 82){\line(-1,2){11.7}}
\put(121, 98){\circle{4}}
\put(126, 96){\scriptsize{$4_3$}}
\put(121,100){\line(0,1){14}}
\put(121,116){\circle{4}}
\put(126,114){\scriptsize{$9_2$}}
\put( 91,107){\circle{4}}
\put( 79,105){\scriptsize{$8_2$}}
\put(106,132){\line(1,-1){14}}
\put(104,132){\line(-1,-2){11.7}}
\put(105,134){\circle{4}}
\put( 80,132){\scriptsize{\fbox{$16_1$}}}
\put(106,136){\line(1,1){14}}
\put(104,136){\line(-1,2){11.7}}
\put(121,152){\circle{4}}
\put(126,150){\scriptsize{$9_3$}}
\put(121,154){\line(0,1){14}}
\put(121,170){\circle{4}}
\put(126,168){\scriptsize{$4_4$}}
\put( 91,161){\circle{4}}
\put( 79,159){\scriptsize{$8_1$}}
\put(106,186){\line(1,-1){14}}
\put(104,186){\line(-1,-2){11.7}}
\put(105,188){\circle{4}}
\put( 84,186){\scriptsize{\fbox{$1_3$}}}
\put(105,190){\line(0,1){16}}
\put(105,208){\circle{4}}
\put(110,206){\scriptsize{$4_2$}}
\put(105,210){\line(0,1){16}}
\put(105,228){\circle{4}}
\put(110,226){\scriptsize{$2_3$}}
\put(105,230){\line(0,1){16}}
\put(105,248){\circle{4}}
\put(110,246){\scriptsize{$1_1$}}
\put(174,  0){$2a>b>a$}
\put(195, 20){\circle{4}}
\put(200, 18){\scriptsize{$1_4$}}
\put(195, 22){\line(0,1){10}}
\put(195, 34){\circle{4}}
\put(200, 32){\scriptsize{$2_4$}}
\put(195, 36){\line(0,1){10}}
\put(195, 48){\circle{4}}
\put(200, 46){\scriptsize{$4_5$}}
\put(195, 50){\line(0,1){10}}
\put(195, 62){\circle{4}}
\put(200, 60){\scriptsize{$2_2$}}
\put(195, 64){\line(0,1){10}}
\put(195, 76){\circle{4}}
\put(181, 74){\scriptsize{$9_4$}}
\put(197, 77){\line(3,2){15.7}}
\put(215, 88){\circle{4}}
\put(220, 86){\scriptsize{$8_4$}}
\put(215, 90){\line(0,1){10}}
\put(215,102){\circle{4}}
\put(220,100){\scriptsize{$1_2$}}
\put(215,104){\line(0,1){10}}
\put(215,116){\circle{4}}
\put(220,114){\scriptsize{$4_3$}}
\put(215,118){\line(0,1){10}}
\put(215,130){\circle{4}}
\put(220,128){\scriptsize{$9_2$}}
\put(195,142){\circle{4}}
\put(168,140){\scriptsize{\fbox{$16_1$}}}
\put(197,141){\line(3,-2){15.7}}
\put(180,109){\circle{4}}
\put(167,107){\scriptsize{$8_2$}}
\put(180,111){\line(1,2){14.5}}
\put(180,107){\line(1,-2){14.5}}
\put(197,143){\line(3,2){15.7}}
\put(215,154){\circle{4}}
\put(220,152){\scriptsize{$9_3$}}
\put(215,156){\line(0,1){10}}
\put(215,168){\circle{4}}
\put(220,166){\scriptsize{$4_4$}}
\put(215,170){\line(0,1){10}}
\put(215,182){\circle{4}}
\put(220,180){\scriptsize{$1_3$}}
\put(215,184){\line(0,1){10}}
\put(215,196){\circle{4}}
\put(220,194){\scriptsize{$8_3$}}
\put(195,208){\circle{4}}
\put(183,206){\scriptsize{$9_1$}}
\put(197,207){\line(3,-2){15.7}}
\put(180,175){\circle{4}}
\put(168,173){\scriptsize{$8_1$}}
\put(180,177){\line(1,2){14.5}}
\put(180,173){\line(1,-2){14.5}}
\put(195,210){\line(0,1){10}}
\put(195,222){\circle{4}}
\put(200,220){\scriptsize{$2_1$}}
\put(195,224){\line(0,1){10}}
\put(195,236){\circle{4}}
\put(200,234){\scriptsize{$4_2$}}
\put(195,238){\line(0,1){10}}
\put(195,250){\circle{4}}
\put(200,248){\scriptsize{$2_3$}}
\put(195,252){\line(0,1){10}}
\put(195,264){\circle{4}}
\put(200,262){\scriptsize{$1_1$}}
\put(284,  0){$b>2a$}
\put(295, 20){\circle{4}}
\put(300, 18){\scriptsize{$1_4$}}
\put(295, 22){\line(0,1){16}}
\put(295, 40){\circle{4}}
\put(300, 38){\scriptsize{$2_4$}}
\put(296, 42){\line(1,1){15.5}}
\put(294, 42){\line(-1,1){15.5}}
\put(277, 59){\circle{4}}
\put(264, 57){\scriptsize{$1_2$}}
\put(313, 59){\circle{4}}
\put(318, 57){\scriptsize{$4_5$}}
\put(296, 76){\line(1,-1){15.5}}
\put(294, 76){\line(-1,-1){15.5}}
\put(295, 78){\circle{4}}
\put(301, 75){\scriptsize{$8_4$}}
\put(296, 80){\line(3,2){15.5}}
\put(294, 80){\line(-1,1){17}}
\put(313, 92){\circle{4}}
\put(318, 90){\scriptsize{$9_4$}}
\put(311, 93){\line(-5,4){33.5}}
\put(313, 94){\line(0,1){14}}
\put(313,110){\circle{4}}
\put(318,108){\scriptsize{$2_2$}}
\put(313,112){\line(0,1){14}}
\put(313,128){\circle{4}}
\put(318,126){\scriptsize{$8_2$}}
\put(296,140){\line(3,-2){15.5}}
\put(295,142){\circle{4}}
\put(266,140){\scriptsize{\fbox{$16_1$}}}
\put(276, 99){\circle{4}}
\put(264, 97){\scriptsize{$4_3$}}
\put(276,101){\line(0,1){18}}
\put(276,121){\circle{4}}
\put(264,119){\scriptsize{$9_2$}}
\put(294,140){\line(-1,-1){17}}
\put(296,144){\line(3,2){15.5}}
\put(294,144){\line(-1,1){17}}
\put(313,156){\circle{4}}
\put(318,154){\scriptsize{$8_1$}}
\put(313,158){\line(0,1){14}}
\put(313,174){\circle{4}}
\put(318,172){\scriptsize{$2_1$}}
\put(313,176){\line(0,1){14}}
\put(313,192){\circle{4}}
\put(318,190){\scriptsize{$9_1$}}
\put(311,191){\line(-5,-4){33.5}}
\put(296,204){\line(3,-2){15.5}}
\put(295,206){\circle{4}}
\put(301,204){\scriptsize{$8_3$}}
\put(276,163){\circle{4}}
\put(264,161){\scriptsize{$9_3$}}
\put(276,165){\line(0,1){18}}
\put(276,185){\circle{4}}
\put(264,183){\scriptsize{$4_4$}}
\put(294,204){\line(-1,-1){17}}
\put(296,208){\line(1,1){15.5}}
\put(294,208){\line(-1,1){15.5}}
\put(277,225){\circle{4}}
\put(264,223){\scriptsize{$1_3$}}
\put(313,225){\circle{4}}
\put(318,223){\scriptsize{$4_2$}}
\put(296,242){\line(1,-1){15.5}}
\put(294,242){\line(-1,-1){15.5}}
\put(295,244){\circle{4}}
\put(300,242){\scriptsize{$2_3$}}
\put(295,246){\line(0,1){16}}
\put(295,264){\circle{4}}
\put(300,262){\scriptsize{$1_1$}}
\end{picture}
\end{center}
{\small A box indicates a two-sided cell with several constructible 
representations, see Table~\ref{lrgraph1}. Otherwise, the two-sided cell 
has only one irreducible, constructible respresentation.}
\end{table}

\begin{table}[htbp] \caption{Left cell representations in type $F_4$} 
\label{lrgraph1}
{\small \begin{center}
$\renewcommand{\arraystretch}{1.3}
\begin{array}{l@{\hspace{1mm}}l} \hline & a=b \\ \hline
\mbox{\fbox{$4_2$}}:  
&2_3 {+} 4_2 ,\\
&2_1 {+} 4_2  \\
\mbox{\fbox{$12_1$}}: &
9_3 {+} 6_1 {+} 12_1 {+} 4_4 {+} 16_1, \\
&9_2 {+} 6_1 {+} 12_1 {+} 4_3 {+} 16_1,  \\
&4_1 {+} 9_2 {+} 9_3 {+} 6_2 {+} 12_1 {+} 2{\cdot} 16_1,  \\
&1_3 {+} 2 {\cdot} 9_3 {+} 6_2 {+} 12_1 {+} 4_4 {+} 16_1,  \\
&1_2 {+} 2 {\cdot} 9_2 {+} 6_2 {+} 12_1 {+} 4_3 {+} 16_1, \\
\mbox{\fbox{$4_4$}}: 
&2_4 {+} 4_5,\\
&2_2 {+} 4_5   \\
\hline \end{array}\quad
\begin{array}{l@{\hspace{1mm}}l} \hline & b=2a\\\hline
\mbox{\fbox{$1_3$}}: 
& 1_3 {+} 8_3,\\
&2_1 {+} 9_1,\\
&9_1 {+} 8_3\\
\mbox{\fbox{$16_1$}}:
&6_1 {+} 12_1 {+} 16_1\\
&6_2 {+} 12_1 {+} 16_1\\
&4_1 {+} 16_1\\
\mbox{\fbox{$1_2$}}: 
&1_2 {+} 8_4,\\
&2_2 {+} 9_4,\\
&9_4 {+} 8_4\\
\hline \end{array}\quad
 \begin{array}{l@{\hspace{1mm}}l} \hline & b\not\in \{a,2a\}\\ \hline
\mbox{\fbox{$16_1$}}:
&6_1 {+} 12_1 {+} 16_1,\\
&6_2 {+} 12_1 {+} 16_1,\\
&4_1 {+} 16_1\\
\hline\\\\\\\\\\\\\end{array}$
\end{center}}
\end{table}

Thus, we have finally covered all cases of unequal parameters. A detailed
analysis of the partition of left cells obtained in each case leads us 
to the following result.

\begin{cor} \label{corconf4} Let $L=L_{a,b}$ and $L'=L_{a',b'}$ be two 
weight functions on $W$ such that $b\geq a>0$ and $b'\geq a'>0$.
Then $L,L'$ define the same partition of $W$ into left cells if and only
if $L,L'\in \cL_i$ for $i\in \{0,1,2,3\}$, where $\cL_i$ are defined
as follows:
\begin{alignat*}{3}
\cL_0 &= \{ (c,c,c,c) && \mid c>0\}, \\
\cL_1 &= \{ (c,c,2c,2c) &&\mid c>0\},\\
\cL_2 &= \{ (c,c,d,d) &&\mid 2c>d>c>0\},\\
\cL_3 &= \{ (c,c,d,d) &&\mid d>2c>0\}.
\end{alignat*}
The partial order relation $\leq_{LR}$ on two-sided cells and the
left cell representations are given in Tables~\ref{lrgraph} and 
\ref{lrgraph1}. In all cases, the left cell representations are 
precisely the constructible representations, as defined in 
\cite[Chap.~22]{Lusztig03}. Furthermore, the statements in 
Theorems~\ref{typf4b} and \ref{typf4c} hold for any weight function~$L$.
\end{cor}

Note that the list of constructible representations given in 
\cite{Lusztig03}, \S 22.27, Case~1, has to be corrected as specified in
Table~\ref{lrgraph1}; see Remark~\ref{rem2} below.

\begin{proof} Let $L=L_{c,d}$ be any weight function on $W$ where
$d\geq c>0$. In addition to the results obtained in 
Lemmas~\ref{lem41}--\ref{lem47}, we use our {\sf CHEVIE} program to 
compute all the required data in the cases where 
\begin{equation*}
 \{c,d\}\in \{(1,4),(1,3),(2,5),(1,2),(2,3),(3,4),(1,1)\}.\tag{1}
\end{equation*}
Then, by Remark~\ref{rem1}, we have  covered all equivalences classes
of weight functions on $W$. In each of the above cases, our {\sf CHEVIE} 
program has automatically computed the preorder relations $\leq_L$ and 
$\leq_{LR}$ and checked that Theorem~\ref{typf4c} holds. Furthermore, by 
inspection of the partitions into left cells obtained in the various cases, 
we find the above four classes of weight functions $\cL_i$ ($0\leq i\leq 3$).
The decompositions of the left cell representations are determined by 
explicit computations using the character table of $W$. By inspection, 
we see that the left cell representations are precisely the constructible 
representations as determined by Lusztig \cite[\S 22.27]{Lusztig03}
(modulo the error in Case~1 in Lusztig's list). 

It remains to prove the statements in Theorem~\ref{typf4b}, concerning
the distinguished involutions. For this purpose, we use a similar procedure
as before, beginning with a total order $\Gamma_+ \subset \Gamma$ as
specified in Lemma~\ref{lem41}. But now we have to work with the larger
set $\Gamma_+'(W)$ defined in Remark~\ref{distinv} in each step and make
sure that ($*^\prime$) holds. For example, the analogue of Lemma~\ref{lem41}
now reads: 

{\em Let $\Gamma_+\subset \Gamma$ be a pure lexicographic order as in 
Lemma~\ref{lem41}. Then condition ($*^\prime$) in Remark~\ref{distinv} 
holds provided that $b/a>9$.}

Then we continue with an anologue of Lemma~\ref{lem42} and so on. Thus, 
there will be more cases to be considered, but the whole argument is 
basically the same. We omit the details. Once this is done, one can argue 
as follows. Let $C$ be a left cell of $W$ (with respect to a total order 
$\Gamma_+ \subset \Gamma$ similar to one of the cases in 
Lemmas~\ref{lem41}--\ref{lem47}). By inspection, one checks that the 
following holds:
\begin{center}
{\em There exists a (unique) $d_0\in C$ such that $\delta_{d_0}^{-1}
\delta_w \in \Gamma_+$ for every $w\in C\setminus\{d_0\}$}.  
\end{center}
(Here, $\delta_w$ is defined as in Remark~\ref{distinv}.) Thus, we may
regard $d_0$ as a {\em distinguished involution} in $C$. Now, the fact that
condition ($*^\prime$) in Remark~\ref{distinv} holds in these cases shows 
that the function $w\mapsto \Delta(w)$ restricted to $C$ also reaches its 
minimum at $d_0\in C$ and that $\Delta(w)>\Delta(d_0)$ for all $w\in 
C\setminus \{d_0\}$.
\end{proof}

\begin{rem} \label{rem0} By inspection of the results obtained in 
Corollary~\ref{corconf4} and its proof, we find the following:
\begin{itemize}
\item Let $C$ be a left cell with respect to a weight function in $\cL_0$. 
Then $C$ is a union of left cells with respect to a weight function in 
$\cL_3$. 
\item Let $C$ be a left cell with respect to a weight function in $\cL_1$. 
Then $C$ is a union of left cells with respect to a weight function in
$\cL_3$, and $C$ also is a union of left cells with respect to a weight 
function in $\cL_4$.
\end{itemize}
Such a behaviour has been conjectured by Bonnaf\'e to hold in general
(private communication). For type $I_2(m)$, this is easily verified 
using the results in \cite[Chap.~8]{Lusztig03}. One can also check that
this holds for type $B_3$ and $B_4$.
\end{rem}

\begin{rem} \label{rem2} Consider the case where $b=2a>0$. Lusztig states 
in \cite[22.27, Case~1]{Lusztig03} that $1_3\oplus 2_1$ and $1_2\oplus 2_2$
are constructible. However, these representations are not constructible. 
(In fact, we just have to omit them from the list given by Lusztig.)
Let us add some details about this. The $\mathbf{a}$-invariants of the 
irreducible representations of $W$ are given by:
\[\begin{array}{ccccccccccccccc} \hline \mathbf{a} & 
 0& a& 2a& 3a& 5a& 6a& 7a& 10a& 11a& 12a& 15a& 20a& 25a& 36a \\ \hline \rho 
 &1_1&2_3&4_2&1_3& 4_4 & 9_3& 4_1 &9_2&4_3&8_2&1_2& 4_5 & 2_4 & 1_4 \\
 &&&&2_1 && 8_1 & 6_1 &  &&&2_2& & & \\
 &&&&9_1 && & 6_2 &  &&&9_4& & & \\
 &&&&8_3 && & 12_1 &  &&&8_4& & & \\
 &&&& && & 16_1 &  &&&& & & \\\hline
\end{array}\]
For $i\in \{1,2,3,4\}$, let $W_i$ be the parabolic subgroup of $W$
generated by $S\setminus\{s_i\}$.  The maximal $\mathbf{a}$-invariant
of a representation of $W_i$ (for $i=1,2,3,4$) is given by $15a$, $7a$, $6a$
or $12a$, respectively. Furthermore, that maximal value is reached only
at the sign representation. Thus, since the restriction of $1_2$ to 
$W_i$ is not the sign representation, we conclude that $1_2$ cannot occur
in the $J$-induction of any representation of any $W_i$. Hence 
$1_3$ (obtained from $1_2$ by tensoring with sign) must occur in the 
$J$-induction from some proper parabolic subgroup. Now, the restriction 
of $1_3$ to $W_1$ (type $C_3$) is given by $(\varnothing,3)$. Furthermore, 
this representation is constructible.
The restriction of $1_3$ to $W_2$ (type $A_1 \times A_2$) is given 
by $(11) \boxtimes (3)$. Furthermore, this representation is constructible.
The restriction of $1_3$ to $W_3$ (type $A_2 \times A_1$) is given by 
$(111) \boxtimes (2)$. Furthermore, this representation is constructible.
The restriction of $1_3$ to $W_4$ (type $B_3$) is given by $(111,
\varnothing)$. Furthermore, the representation $(111,\varnothing)+(11,1)$ 
is constructible, and this is the only constructible representation in 
which $(111,\varnothing)$ occurs; see \cite[Chap.~22]{Lusztig03}. We have
\begin{gather*}
\mbox{J}_{W_1}^W\bigl((\varnothing,3)\bigr)=2_3,\qquad
\mbox{J}_{W_3}^W\bigl((111) \boxtimes (2)\bigr)=1_3 \oplus 8_3,\\
\mbox{J}_{W_2}^W\bigl((11) \boxtimes (3)\bigr)= 2_3,\qquad 
\mbox{J}_{W_4}^W\bigl((111,\varnothing)+(11,1)\bigr)= 1_3 \oplus 8_3.
\end{gather*}
Thus, $1_3\oplus 8_3$ is the only constructible representation of $W$
in which $1_3$ occurs.
\end{rem}

\begin{rem} \label{rem3}
The case $b=2a$ in type $F_4$ also shows that, in general, there no longer 
exist representations which would have similar properties as the ``special'' 
representations in the equal parameter case (see \cite[\S 12]{Lu82}). Indeed,
consider the two-sided cell containing $1_3$. Then the three constructible
representations belonging to that two-sided cell do not have an irreducible
constituent in common.
\end{rem}


\end{document}